\documentclass{svmult-ddm}

\usepackage{mathptmx}
\usepackage{helvet}
\usepackage{courier}
\usepackage{type1cm}

\usepackage{graphicx}

\usepackage[bottom]{footmisc}

\usepackage{amssymb}
\usepackage{amsfonts}
\usepackage{amsbsy}
\usepackage{amscd}
\usepackage{amstext}
\usepackage{dsfont}
\usepackage[english]{babel}
\usepackage{color}
\usepackage{graphics}
\usepackage{epsfig}
\usepackage{subfigure}
\usepackage{wrapfig}
\usepackage{psfrag}
\usepackage{color}
\usepackage{url}
\usepackage{verbatim}
\usepackage{algorithm}
\usepackage{algorithmic}
 
\newcommand{\R}{\mathbb{R}}
\newcommand{\N}{\mathbb{N}}

\newcommand{\Q}{\mathbb{Q}}
\newcommand{\cR}{\mathcal{R}}
\newcommand{\cI}{\mathcal{I}}
\newcommand{\cG}{\mathcal{G}}
\newcommand{\cF}{\mathcal{F}}
\newcommand{\cT}{\mathcal{T}}
\newcommand{\cO}{\mathcal{O}}

\newcommand{\average}[1]{\ensuremath{\left\{#1\right\}}}
\newcommand{\jump}[1]{\ensuremath{\left[#1\right]}}
\newcommand{\skp}[1]{\ensuremath{\left\langle{#1}\right\rangle}}

\newcommand{\discrete}[1]{{\bf{#1}}}

\newcommand{\BIGOP}[1]{\mathop{\mathchoice%
{\raise-0.22em\hbox{\huge $#1$}}%
{\raise-0.05em\hbox{\Large $#1$}}{\hbox{\large $#1$}}{#1}}}
\newcommand{\bigtimes}{\BIGOP{\times}}

\begin{document}

\title*{Robust Preconditioners for DG-Discretizations with Arbitrary Polynomial Degrees}
\titlerunning{Robust Preconditioners for Spectral DG}
\author{Kolja Brix\inst{1} \and Claudio Canuto\inst{2} \and Wolfgang Dahmen\inst{1}}
\institute{Institut f\"ur Geometrie und Praktische Mathematik, RWTH Aachen, Templergraben 55, 52056 Aachen, Germany
\texttt{\{brix,dahmen\}@igpm.rwth-aachen.de}
\and Dipartimento di Scienze Matematiche, Politecnico di Torino, Corso Duca degli Abruzzi 24, 10129 Torino, Italy,
\texttt{claudio.canuto@polito.it}
}

\maketitle

\abstract*{Discontinuous Galerkin (DG) methods offer an enormous flexibility regarding local grid refinement and variation of polynomial degrees for a variety of different problem classes. With a focus on diffusion problems, we consider DG discretizations for elliptic boundary value problems, in particular the efficient solution of the linear systems of equations that arise from the Symmetric Interior Penalty DG method. We announce a multi-stage preconditioner which produces uniformly bounded condition numbers and aims at supporting the full flexibility of DG methods under mild grading conditions. The constructions and proofs are detailed in an upcoming series of papers by the authors. Our preconditioner is based on the concept of the auxiliary space method and techniques from spectral element methods such as Legendre-Gau{\ss}-Lobatto grids. The presentation for the case of geometrically conforming meshes is complemented by numerical studies that shed some light on constants arising in four basic estimates used in the second stage.
}


\section{Introduction}\label{sec:Intro}
\vspace*{-2mm}
Discontinuous Galerkin (DG) methods offer an enormous flexibility regarding local grid refinement and variation of polynomial degrees rendering such concepts powerful discretization tools which have proven to be well-suited for a variety of different problem classes. While initially the main focus has been on transport problems like hyperbolic conservation laws, interest has meanwhile shifted towards diffusion problems. Specifically, we focus here on the efficient solution of the linear systems of equations that arise from the Symmetric Interior Penalty DG method applied to elliptic boundary value problems.~\cite{Ayuso2012} The principal objective is to develop robust preconditioners for the full ``DG-flexibility'' which means to obtain uniformly bounded condition numbers for locally refined meshes and arbitrarily (subject to mild grading conditions) varying polynomial degrees at the expense of linearly scaling computational work. A first step towards that goal has been made in \cite{BCC+2012} treating the case of geometrically conforming meshes but arbitrarily large variable polynomial degrees which already exposes major principal obstructions. In this paper we complement this work by detailed studies of several issues arising in \cite{BCC+2012}.

To our knowledge the only concept yielding full robustness with respect to polynomial degrees is based on {\em Legendre-Gau{\ss}-Lobatto} (LGL) quadrature. Specifically, in the framework of {\em auxiliary space methods} low order finite element discretizations on LGL-grids can be used to precondition high order polynomial discretizations. However, when dealing with variable degrees the possible non-matching of such grids at element interfaces turns out to severely obstruct in general the design of efficient preconditioners. To overcome these difficulties we propose in \cite{BCC+2012} a concatenation of auxiliary space preconditioners. In the first stage the spectral DG formulation (\text{\bf SE-DG}) is transferred to a spectral continuous Galerkin formulation (\text{\bf SE-CG}). In the second stage we proceed from here to a finite element formulation on a specific dyadic grid (\text{\bf DFE-CG}) which is associated with an LGL-grid but belongs to a nested hierarchy. The latter problem can then be tackled by a multilevel wavelet preconditioner presented in forthcoming work. The overall path of our iterated auxiliary space preconditioner therfore is $\text{\bf SE-DG} \to \text{\bf SE-CG} \to \text{\bf DFE-CG}$. It should be noted that a natural alternative is to combine the first stage with a domain decomposition substructuring preconditioner as proposed in \cite{CPP2012} admitting a mild growth of condition numbers with respect to the polynomial degree.

We are content here for most part of the paper with brief pointers to the detailed analysis in \cite{BCC+2012}, \cite{BCC+2012a} and~\cite{Brix2013} to an extent needed for the present discussion.

Section~\ref{sec:ModelProblem} introduces our model problem, the LGL technique is explained in Section~\ref{sec:LGL}. The auxiliary space method is detailed in Section~\ref{sec:ASM}, while Sections~\ref{sec:Stage1} and~\ref{sec:Stage2} consider stages 1 and 2 of our preconditioner. Finally in Section~\ref{sec:NumericalExperiments} we give some numerical experiments that shed light on the constants that arise in four basic inequalities used in the second stage.
\vspace*{-6mm}


\section{Model problem and Discontinuous Galerkin formulation}\label{sec:ModelProblem}
\vspace*{-2mm}
Given a bounded Lipschitz domain $\Omega \subset \R^d$ with piecewise smooth boundary we consider as a simple model problem the weak formulation: find $u \in H^{1}_{0}(\Omega)$ such that
\begin{eqnarray*}
a(u,v):=\int_{\Omega} \nabla u \cdot \nabla v \ \mathrm{d}x=\skp{f,v}, \quad v \in H^{1}_{0}(\Omega)
\end{eqnarray*}
of Poisson's equation $-\Delta u=f$ on $\Omega$ with zero Dirichlet boundary conditions $u=0$ on $\partial \Omega$. For simplicity, we assume that $\bar \Omega$ is the union of a collection $\cR$ of finitely many (hyper-)rectangles, which at most overlap with their boundaries. More complex geometries can be handled by isoparametric mappings. By $\cF_l(R)$ we denote the $l$-dimensional facets of a (hyper-)rectangle $R$ and by $\cF_l=\cup_{R \in \cR} \cF_l(R)$ the union of all these objects. Let $H_k(R)$ be the side length of $R$ in the $k$-th coordinate direction.

The polynomial degrees used in each cell $R$ are defined as $p=(p_k)_{k=1}^d$, where $p_k$ is the polynomial degree in the $k$-th coordinate direction. We introduce the piecewise constant function $\delta=(H,p)$ that collects the $hp$ approximation parameters. On $\delta$ we impose mild grading conditions, see~\cite{BCC+2012} for the details.

For the spectral discretization of our model problem, we choose the DG spectral ansatz space $V_{\delta}:= \left\{ v \in L^2(\Omega) : v|_R \in \Q_p(R) \ \text{for all} \ R \in \cR \right\}$, where $\Q_p(R)$ is the tensor space of all polynomials of degree at most $p$ on the (hyper-)rectangle $R$.

We employ the standard notation of DG methods for jumps and averages on the mesh skeleton and on $\partial\Omega$.

The {\em Symmetric Interior Penalty Discontinuous Galerkin} method (SIPG) $a_{\delta}(u,v)=\skp{f,v}$ for all $v \in V_{\delta}$ with the SIPG bilinear form
\begin{eqnarray*}
a_{\delta}(u_{\delta},v_{\delta}):=
  \sum_{R \in \cR} (\nabla u_{\delta}, \nabla v_{\delta})_{R}
+ \sum_{F \in \cF} ( -(\average{\nabla u_{\delta}}, \jump{v_{\delta}})_{F}
                     -(\jump{u_{\delta}}, \average{\nabla v_{\delta}})_{F} )\\
+ \sum_{F \in \cF} \gamma \omega_F (\jump{u_{\delta}}, \jump{v_{\delta}})_{F}
=(f,v_{\delta})_{\Omega}, \quad v_{\delta} \in V_{\delta}
\end{eqnarray*}
with $\omega_F := \max \left\{ \omega_{F,R^{-}}, \omega_{F,R^{+}} \right\}$ for internal faces $F$ and $\omega_{F,R^{\pm}} := \frac{p_k(R^{\pm}) (p_k(R^{\pm})+1)}{H_k(R^{\pm})}$. For boundary faces $F \subset \partial \Omega $ we set $\omega_{F,R} := \frac{p_k(R) (p_k(R)+1)}{H_k(R)}$.


\section{Legendre-Gau{\ss}-Lobatto (LGL) grids}\label{sec:LGL}
\vspace*{-2mm}
Denoting by $(\xi_i)_{i=1}^{p-1}$ the zeros of the first derivative of the $p$-th Legendre polynomial $L_p$, (in ascending order), and setting $\xi_0=-1$ and $\xi_p=1$, $\cG_p=(\xi_i)_{0 \le i \le p}$ is the Legendre-Gau{\ss}-Lobatto (LGL) grid of degree $p$ on the reference interval $\hat I=[-1,1]$, see e.g. \cite{CHQ+2006}. In combination with appropriate LGL weights $(w_i)_{0 \le i \le p}$ the LGL points of order $p$ can be interpreted as quadrature points of a quadrature rule of exactness $2p-1$. In~\cite{BCC+2012a} we prove quasi-uniformity of the LGL-grids $(\cG_p)_{p \in \N}$, i.e., $\frac{h_{i+1},p}{h_{i},p}$ remains uniformly bounded independent of $p$, where $h_{i}=|\xi_{i}-\xi_{i-1}|$ for $1 \le i \le p-1$.

The particular relevance of tensor product LGL-grids for preconditioners for spectral element discretizations lies in the two norm equivalences (see ~\cite{CHQ+2006})
\begin{eqnarray}
\big\|\varphi\big\|_{H^i(R)} \eqsim \big\|\cI_{h,p}^{R} \varphi\big\|_{H^i(R)} \quad \text{for all} \quad \varphi \in \Q_{p}(R), \quad i \in \{0,1\},
\end{eqnarray}
which hold uniformly for any $d$-dimensional hypercube $R=\bigtimes_{k=1}^d I_k$ where $\cI^R_{h,p}$ is the piecewise multi-linear interpolant on the tensor product LGL-grid.
\vspace*{-4mm}


\section{Abstract theory: Auxiliary Space Method}\label{sec:ASM}
\vspace*{-2mm}
The auxiliary space method (ASM)~\cite{Oswald1996,Xu1996,Xu1992} is a powerful concept for the construction of preconditioners that can be derived from the {\em fictitious space lemma} \cite{Nepomnyaschikh1992,Nepomnyaschikh1990,Oswald1996}.

Given a problem $a(u,v)=f(v)$ for all $v \in V$ on the linear space $V$ equipped with a bilinear form $a(\cdot,\cdot): V \times V \rightarrow \R$, we seek an {\em auxiliary space} $\tilde V$ with an {\em auxiliary form} $\tilde a(\cdot,\cdot): \tilde V \times \tilde V \rightarrow \R$ that is in some sense close to the original one but easier to solve. Note that we neither require $V \subset \tilde V$ nor $\tilde V \subset V$ which is important in the context of non-conforming discretizations. Therefore on the sum $\hat V = V + \tilde V$ we need in general another version $\hat a(\cdot,\cdot): \hat V \times \hat V \rightarrow \R$ as well as a second form $b(\cdot,\cdot): \hat V \times \hat V \rightarrow \R$ which dominates $a$ on $V$ and plays the role of a smoother. The required closeness of the spaces $V$ and $\tilde V$ is described with the aid of two linear operators $Q: \tilde V \rightarrow V$ and $\tilde Q: V \rightarrow \tilde V$. Specifically, these operators have to satisfy certain direct estimates involving the above bilinear forms. For the details on the ASM conditions see~\cite{Oswald1996}.

\begin{lemma}[Stable Splitting \cite{Oswald1996}]
Under the assumptions of the ASM, we have the following stable splitting
\begin{eqnarray*}
a(v,v) \sim \inf_{w \in V, \tilde v \in \tilde V: \ v=w+Q\tilde v} \left( b(w,w) + \tilde a(\tilde v,\tilde v) \right) \quad \text{for all} \ v \in V.
\end{eqnarray*}
\end{lemma}

The main result of the ASM is given in the following theorem \cite{Oswald1996}.
\begin{theorem}[Auxiliary Space Method]
Let $\discrete{C_{B}}$ and $\discrete{C_{\tilde A}}$ be symmetric preconditioners for $\discrete{B}$ and $\discrete{\tilde A}$, respectively. Let $\discrete{S}$ be the representation of $Q:\tilde V \rightarrow V$. Then $\discrete{C_{A}}:=\discrete{C_{B}}+\discrete{S} \discrete{C_{\tilde A}} \discrete{S}^T$ is a symmetric preconditioner for $\discrete{A}$. Moreover, there exists a uniform constant $C$ such that the spectral condition number of $\discrete{C_{A}} \discrete{A}$ satisfies
\begin{eqnarray*}
\kappa(\discrete{C_{A}} \discrete{A}) \le C \kappa(\discrete{C_{B}}\discrete{B}) \kappa(\discrete{C_{\tilde A}} \discrete{\tilde A}).
\end{eqnarray*}
\end{theorem}

For a given practical application it remains to identify a suitable auxiliary space $\tilde V$, the bilinear forms $\tilde a: \tilde V \times \tilde V \rightarrow \R$ and $\hat a, b: \hat V \times \hat V \rightarrow \R$, as well as the two linear operators $Q$ and $\tilde Q$, such that ASM conditions are satisfied. In addition efficient preconditioners for the ``easier'' auxiliary problems $\discrete{C_{\tilde A}}$ and $\discrete{C_{B}}$ need to be devised. Of course, the rationale is that the complexity to apply $\discrete{C_{\tilde A}}$ and $\discrete{C_{B}}$ should be much lower than solving the original problem.

Note that the operator $\tilde Q$ need {\em not} be implemented but enters only the analysis.


\vspace*{-4mm}
\section{Stage 1: ASM DG-SEM $\rightarrow$ CG-SEM}\label{sec:Stage1}
\vspace*{-2mm}
In the first stage, we choose the largest conforming subspace $\tilde V:=V_\delta \cap H_0^1(\Omega)$ of $V:=V_\delta$ as auxiliary space so that $Q$ can be taken as the canonical injection. The definition of the operator $\tilde Q$ can be found in \cite{BCC+2012}.

The main issue in this stage is the choice of the auxiliary form $b(\cdot,\cdot)$. Using LGL-quadrature combined with an inverse estimate for the partial derivatives in the form $a(\cdot,\cdot)$ we arrive at
\vspace{-3mm}
\begin{eqnarray*}
b(u,v):= \sum_{R \in \cR} \sum_{\xi \in \cG_p(R)} u(\xi) \, v(\xi) \, c_{\xi} W_{\xi},
\quad
W_{\xi}:=\left(\sum_{k=1}^{d} w_{\xi,k}^{-2} \right) w_{\xi,k}.
\end{eqnarray*}
Here the weights $c_{\xi} \sim 1$ are chosen as
\vspace{-2mm}
\begin{eqnarray*}
c_{\xi} := \left\{
\begin{array}{cc}
\beta_1(c_1^2+\gamma \rho_1 \omega_F w_{F,R}/W_{\xi}),& \xi \in \cG_p(F,R), \ F\in \cF_{d-1}(R), \ R \in \cR,\\
\beta_1 c_1^2, & \textnormal{else},
\end{array}\right.
\end{eqnarray*}
where $w_{F,R^{\pm}}$ is the LGL quadrature weight on $F$ seen as a face of $R^{\pm}$ and the parameters $\beta, \rho_1$ can be used to ``tune'' the scheme. The resulting matrix $\discrete{B}$ is diagonal so that the application of $\discrete{C_B}:=\discrete{B}^{-1}$ requires only $\cO(N)$ operations. It is shown in \cite{BCC+2012} that all ASM conditions are satisfied for this choice of $b(\cdot,\cdot)$. Numerical experiments show that the parameters $\beta_1$ and $\rho_1$ can by and large be optimized independently of the polynomial degrees.
\vspace*{-4mm}


\section{Stage 2: CG-SEM $\rightarrow$ CG-DFEM}\label{sec:Stage2}
\vspace*{-2mm}

The second stage involves three major ingredients, namely
\begin{list}{}{%
  \setlength{\labelwidth}{1.7em}
  \setlength{\labelsep}{0.3em}
  \setlength{\leftmargin}{2em}
  \setlength{\rightmargin}{0pt}
  \setlength{\listparindent}{0pt}
  \setlength{\itemindent}{0pt}
  }
\item[(1)] the choice of spaces of piecewise multi-linear finite elements on hierarchies of {\em nested} anisotropic dyadic grids, to permit a subsequent application of efficient multilevel preconditioners,
\item[(2)] the construction of the operators $Q$ and $\tilde Q$, and
\item[(3)] the construction of the auxiliary bilinear form $b(\cdot,\cdot)$.
\end{list}
As for (1), the non-matching of LGL-grids for different degrees $p$ at interfaces prevents us from taking low order finite element spaces as auxiliary space for the high order conforming problem resulting from the first stage. Instead, with each LGL-grid $\cG_p$ we associate a dyadic grid $\cG_{D,p}$, which is roughly generated as follows: starting with the boundary points $\{-1,1\}$ as initial guess we adaptively refine the grid. A subinterval in the grid is bisected into two parts of equal size, if the smallest of the overlapping LGL-subintervals is longer than $\alpha$ times its length. The parameter $\alpha$ therefore controls the mesh size of the dyadic grid. However, for input LGL-grids of different polynomial degrees the resulting dyadic grids are not necessarily nested yet. How to ensure nestedness while keeping the grid size under control is shown in \cite{BCC+2012}. The key quality of the associated dyadic grids $\cG_{D,p}$ is that mutual low order piecewise multi-linear interpolation between the low order finite element spaces on $\cG_p(R), \cG_{D,p}(R)$ is uniformly $H^1$-stable, see \cite{BCC+2012} for the proofs. Denoting by $V_{h,D,p}(R)$ the space of piecewise multi-linear conforming finite elements on $\cG_{D,p}(R)$, we now take $V:=V_\delta \cap H^1_0(\Omega)$ and $\tilde V :=V_{h,D} \cap H^1_0(\Omega)$, where $V_{h,D} = \{{v} \in C^0(\overline{\Omega}) \ : \ \forall R \in \cR\;, \ {v}_{|R}:={v}_R \in V_{h,D,p}(R) \}$.

Concerning (2), the operator $Q$ is defined element-wise as follows. For a given element vertex $z\in \cF_0(R)$ let $p^*$ denote the polynomial degree vector whose $k$th entry is the minimum of the $k$th entries of all degree vectors associated with elements $R'$ sharing $z$ as a vertex. Here a grading of the degrees is important. Let $\Phi_{z} \in \mathbb{Q}_1(R)$ the multi-linear shape function on $R$ satisfying conditions $\Phi_{z}(y)= \delta_{y,z}$ for all $y \in \cF_0(R)$. Then, we define
\begin{eqnarray}
\tilde{v}_z^* := \cI^R_{h,D,p_z^*} \left( \Phi_{z} \tilde{v}_R \right) \in V_{h,D,p_z^*}(R)
\qquad \text{and} \qquad
{v}_z^* = \cI^R_{p_z^*} \, \tilde{v}_z^* \in \mathbb{Q}_{p_z^*}(R) \;,
\end{eqnarray}
where $\cI^R_{h,D,p_z^*}, \cI^R_{p_z^*}$ are the dyadic piecewise multilinear and high order LGL-interpolants on the respective grids. Summing-up over the vertices of $R$, we define
\begin{eqnarray}
\tilde{v}_R^* := \sum_{z \in \cF_0(R)} \tilde{v}_z^* \in V_{h,D,p}(R)
\qquad \text{and} \qquad
Q_R \tilde{v}_R := {v}_R^* := \sum_{z \in \cF_0(R)} {v}_z^* \in \mathbb{Q}_{p}(R) \;.
\end{eqnarray}
The operator $\tilde Q$ is defined analogously with the roles of dyadic and LGL-grids exchanged, see \cite{BCC+2012}.

To finally address (3), for the structure of the form $b(\cdot,\cdot)$ from the first stage the direct estimates in the ASM conditions are no longer valid. It has to be suitably relaxed along the following lines. We make an ansatz of the form
\begin{eqnarray}
\label{eq:b_new}
b(v,w):=\sum_{R\in \cR}\sum_{k=1}^d\Big(\sum_{S_\ell\in \cT_{0,k}(R)}b^0_{R,k,S_\ell}(v,w)
+ \sum_{S_\ell\in \cT_{1,k}(R)}b^1_{R,k,S_\ell}(v,w)\Big),
\end{eqnarray}
where $\cT_{0,k}(R)$ is the collection of those LGL-subcells $S_\ell$, $\ell\in \bigtimes_{k=1}^d \{1, \dots, p_k(R)\} $ with side lengths $h_l^{(\ell_l)}$ in the LGL-grid $\cG_p(R)$ that are {\em strongly anisotropic} according to $ (\max_{l\neq k} h_l^{(\ell_l)})/{h_k^{(\ell_k)}}> C_\textnormal{aspect} $ for a fixed constant $C_\textnormal{aspect}>0$, while $\cT_{1,k}(R)$ is comprised of the remaining ``isotropic'' cells. On the isotropic cells in $\cT_{1,k}(R)$ we use an inverse estimate applied to piecewise multi-linear LGL-interpolants of $v$ and $w$. On the remaining anisotropic cells we retain integrals over the variable involving the partial derivative and use quadrature in the remaining variables. For this auxiliary form $b(\cdot,\cdot)$ and the above operators $Q$ and $\tilde Q$ we can verify all ASM conditions, see \cite{BCC+2012}. Note that the Gramian $\discrete{B}$ is no longer diagonal and we refer to \cite{BCC+2012} for efficient realizations of $\discrete{C_B}$.
\vspace*{-6mm}


\section{Numerical experiments: Constants in the basic interpolation inequalities}\label{sec:NumericalExperiments}
\vspace*{-3mm}
A fundamental role in the proof of the ASM-conditions in the second stage $\text{\bf SE-CG} \to \text{\bf DFE-CG}$ is played by four basic interpolation estimates. In particular, knowing the size of the constants arising in these inequalities and their dependence on the polynomial degrees helps understanding the quantitative effects observed in more complex situations later on.

As before, let $\Phi_{z}$ denote the affine shape function now on the reference interval $\hat I=[-1,1] \subset \R$ satisfying $\Phi_{z}(x)=\delta_{x,z}$ for $x,z \in \{-1,1\}$. By $\cI_{q}$ we denote the polynomial interpolation operator on the LGL-grid $\mathcal{G}_{q}$ for polynomial degree $q$ and by $\cI_{h,D,q}$ the piecewise affine interpolation operator on the dyadic grid $\mathcal{G}_{D,q}$ associated with $\mathcal{G}_{q}$.

A major tool for proving the ASM conditions is given by the following theorem.
\begin{theorem}
Assume that $cp \le q \le p$ for some fixed constant $c>0$. Then we have
\begin{equation}\label{eq:basic0}
| \cI_{q} ( \Phi_{z} v ) |_{H^m(\hat I)} \lesssim \Vert v \Vert_{H^m(\hat I)}
\quad \text{for all} \ v \in \mathbb{Q}_p(\hat I),\ z \in \{-1,1\},\ m \in \{0,1\},
\end{equation}
and
\begin{equation}\label{eq:basic1}
| \cI_{h,D,q} ( \Phi_{z} \tilde v ) |_{H^m(\hat I)} \lesssim \Vert \tilde v
\Vert_{H^m(\hat I)} \quad \text{for all} \ \tilde v \in V_{h,D,p}(\hat I),\ z \in \{-1,1\},\ m \in \{0,1\}.
\end{equation}
\end{theorem}

\begin{figure}
\centering
\subfigure[(\ref{eq:basic0}), $m=0$\label{fig:checkBasicInequalities_3d_2_0}]{\includegraphics[width=0.35\linewidth]{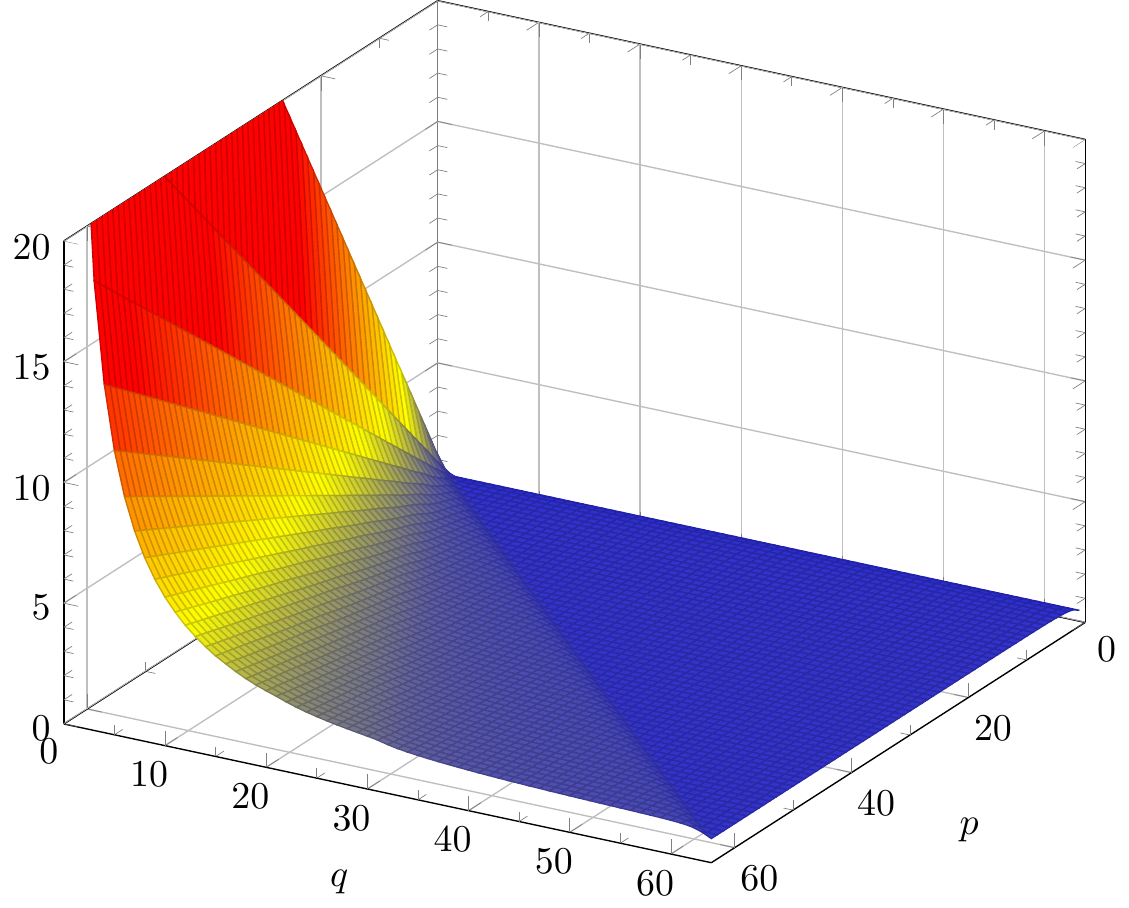}}\hspace{0.1\linewidth}
\subfigure[(\ref{eq:basic0}), $m=1$\label{fig:checkBasicInequalities_3d_2_1}]{\includegraphics[width=0.35\linewidth]{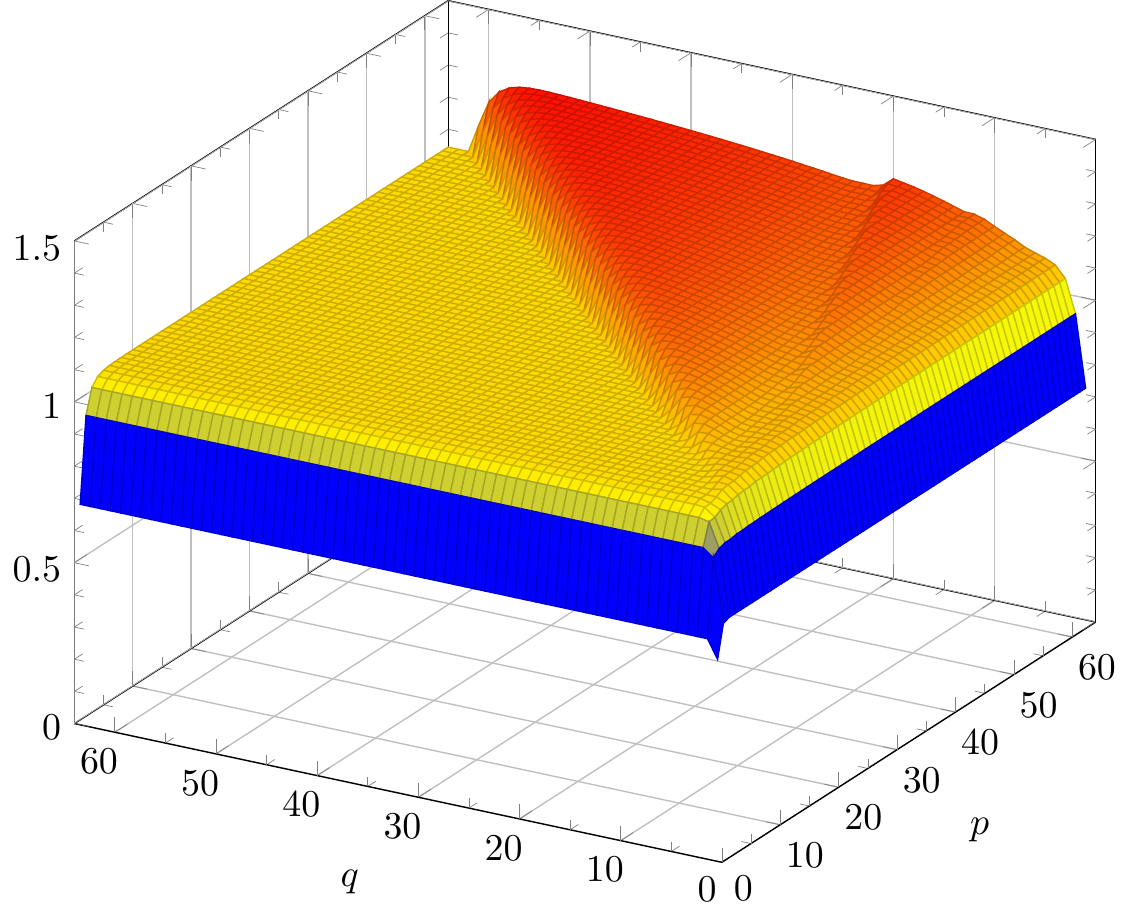}}\\
\subfigure[(\ref{eq:basic1}), $m=0$\label{fig:checkBasicInequalities_3d_1_0}]{\includegraphics[width=0.35\linewidth]{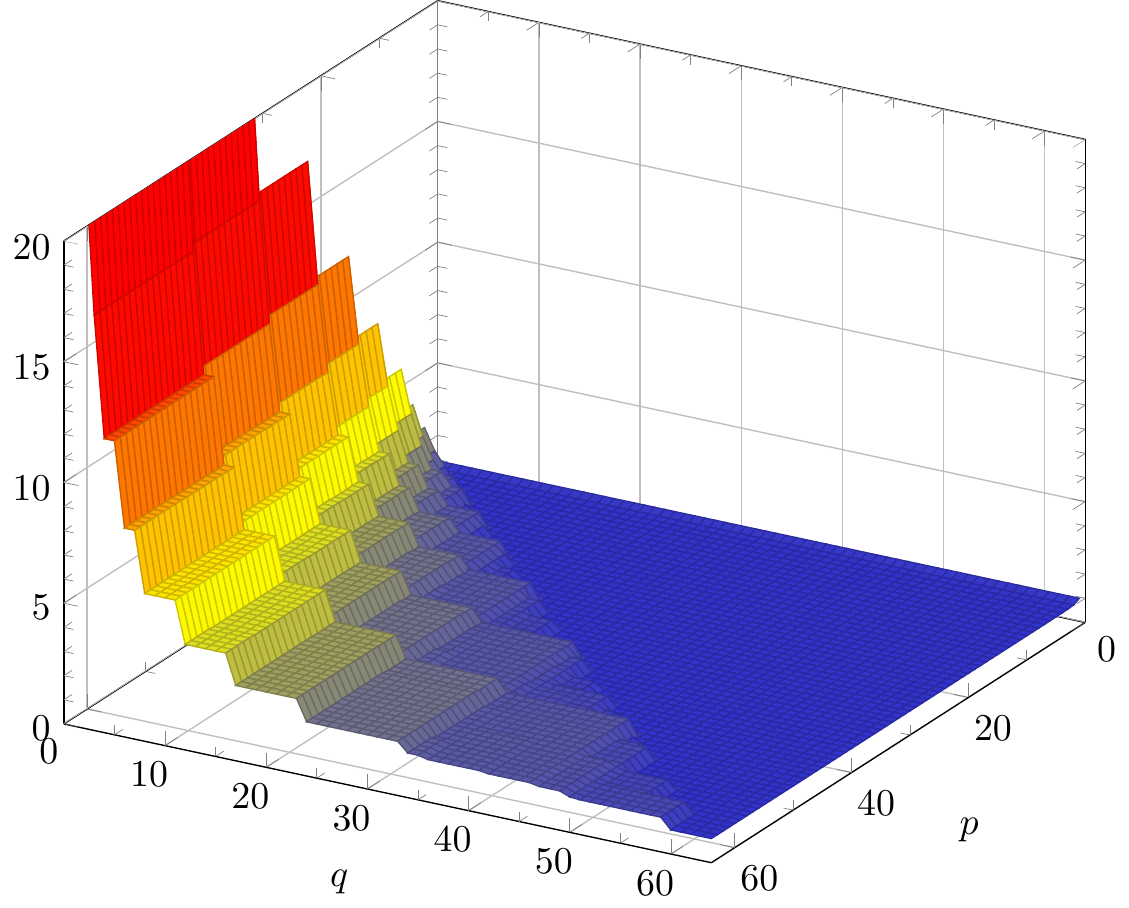}}\hspace{0.1\linewidth}
\subfigure[(\ref{eq:basic1}), $m=1$\label{fig:checkBasicInequalities_3d_1_1}]{\includegraphics[width=0.35\linewidth]{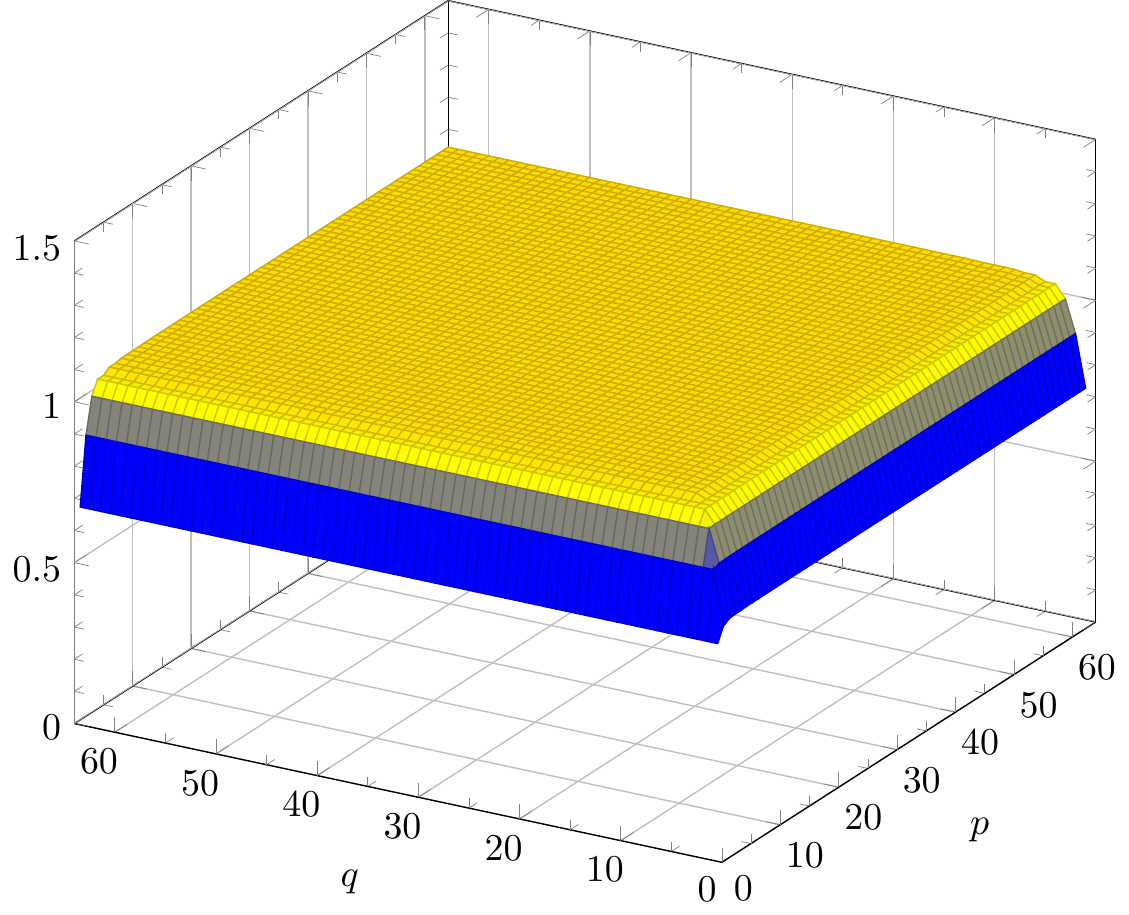}}
\caption{Dependance of the constants in (\ref{eq:basic0}) and (\ref{eq:basic1}) on $p$ and $q$.}
\label{fig:checkBasicInequalities_3d}
\end{figure}

We determine next {\em numerically} the smallest constants that fulfill the inequalities (\ref{eq:basic0}) and (\ref{eq:basic1}). This can be obtained by solving generalized eigenvalue problems for the largest generalized eigenvalue. For all dyadic grids we choose the grid generation parameter $\alpha=1.2$, which balances two effects: on the one hand, the generated auxiliary space is rich enough for a good approximation while on the other hand, to keep the solution of the auxiliary space feasible, the dyadic grid does not have too many degrees of freedom. Figure~\ref{fig:checkBasicInequalities_3d} shows the dependence of the smallest possible constants on the polynomial degrees $p$ and $q$ in the range $1 \le p,q \le 64$.

We observe that the constants in (\ref{eq:basic0}) and (\ref{eq:basic1}) become large for $m=0$ when the quotient $p/q$ increases, but eventually stay bounded as long as $cp \le q \le p$ for a fixed $c>0$. For $m=1$ we find uniform moderate constants in (\ref{eq:basic0}) and (\ref{eq:basic1}) for arbitrary choices of $p$ and $q$. While the nodes in the LGL-grids move gradually with increasing degree the associated dyadic grids change more abruptly which explains the jumps in the graph in Figure~\ref{fig:checkBasicInequalities_3d_1_0}.

\begin{figure}
\centering
\subfigure[m=0]{\includegraphics[width=0.35\linewidth]{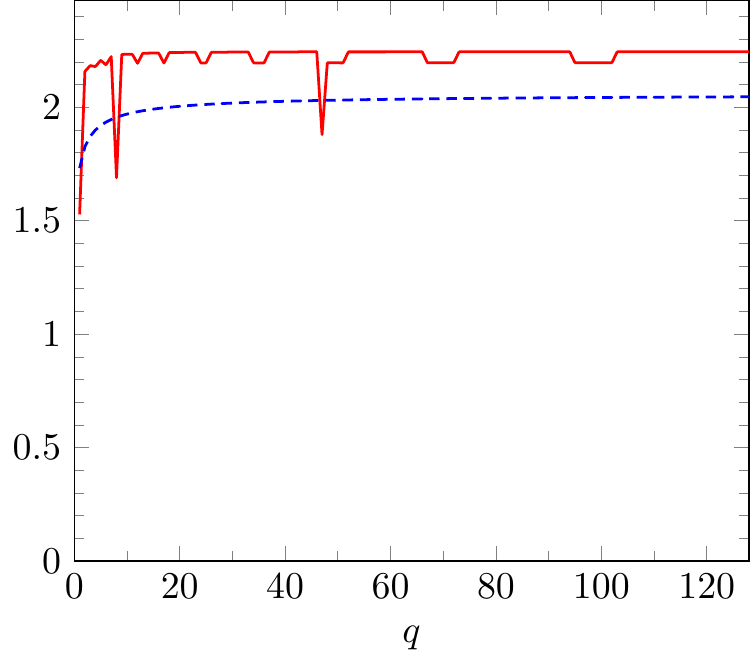}}\hspace{0.1\linewidth}
\subfigure[m=1]{\includegraphics[width=0.35\linewidth]{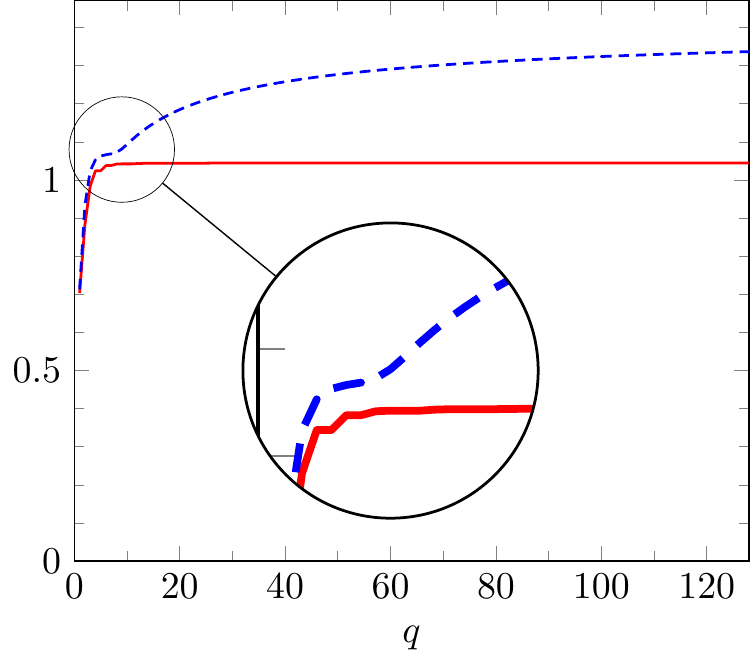}}
\caption{Constants in the basic interpolation inequalities for $p=2q$ (dashed line: (\ref{eq:basic0}), solid line: (\ref{eq:basic1})).}
\label{fig:checkBasicInequalitiesLong}
\end{figure}

We are particularly interested in the behavior of the constants when the quotient of $p$ and $q$ is fixed, i.e., we restrict ourselves to a cross section through the 3-dimensional plots along a line in the $pq$-plane. As an example, we choose $p=2q$ representing strongly varying degrees on adjacent elements. The smallest constants in the inequalities for polynomial degrees $q$ up to $128$ are displayed in Figure~\ref{fig:checkBasicInequalitiesLong}.

While for $m=0$ the constants quickly approach an asymptotic value for both (\ref{eq:basic0}) and for (\ref{eq:basic1}), this is not true for (\ref{eq:basic0}) and $m=1$. In this case we observe a very slow monotonic convergence to its asymptotic limit. Thus for moderate polynomial degrees one still observes a significant growth. Since this estimate is relevant for the ASM conditions on the operator $\tilde Q$ in the second stage, this leads to some growth of the condition number of the preconditioned problem for moderate polynomial degrees and significant inter-element jumps, although it eventually stays uniformly bounded independent of the polynomial degree $q$.
\vspace*{-6mm}


\section{Summary and outlook}
\vspace*{-2mm}
In this paper we sketch a preconditioner for the spectral symmetric interior penalty discontinuous Galerkin method that, under mild grading conditions, is robust in variably arbitrarily large polynomial degrees, announcing detailed results given in \cite{BCC+2012}. The concept is based on the LGL-techniques for spectral methods combined with judiciously chosen nested dyadic grids through an iterated application of the auxiliary space method. A detailed exposition of a multiwavelet preconditioner for the dyadic grid problem, an extension to locally refined grids with hanging nodes, strategies for parallel implementations, and the treatment of jumping coefficients will be presented in forthcoming work.

\vspace*{-3mm}

\begin{acknowledgement}
We thank for the support by the Seed Funds project funded by the Excellence Initiative of the German federal and state governments and by the DFG project 'Optimal preconditioners of spectral Discontinuous Galerkin methods for elliptic boundary value problems' (DA 117/23-1).
\end{acknowledgement}



\begin{thebibliography}{10}
\providecommand{\url}[1]{{#1}}
\providecommand{\urlprefix}{URL }
\expandafter\ifx\csname urlstyle\endcsname\relax
  \providecommand{\doi}[1]{DOI~\discretionary{}{}{}#1}\else
  \providecommand{\doi}{DOI~\discretionary{}{}{}\begingroup
  \urlstyle{rm}\Url}\fi

\bibitem{Ayuso2012}
Ayuso De~Dios, B.: Solvers for discontinuous {G}alerkin methods.
\newblock In: Proceedings of the 21st International Conference on Domain
  Decomposition Methods, Rennes, June 2012 (2012)

\bibitem{Brix2013}
Brix, K.: Robust preconditioners for hp-discontinuous {G}alerkin
  discretizations for elliptic problems.
\newblock Ph.D. thesis, Institut f\"ur Geometrie und Praktische Mathematik,
  RWTH Aachen, in preparation

\bibitem{BCC+2012}
Brix, K., Campos~Pinto, M., Canuto, C., Dahmen, W.: Multilevel preconditioning
  of discontinuous {G}alerkin spectral element methods. {P}art {I}:
  Geometrically conforming meshes.
\newblock IGPM Preprint, RWTH Aachen. In preparation. (2012)

\bibitem{BCC+2012a}
Brix, K., Campos~Pinto, M., Canuto, C., Dahmen, W.: Some properties of
  {L}egendre-{G}au{\ss}-{L}obatto intervals.
\newblock IGPM Preprint, RWTH Aachen. In preparation. (2012)

\bibitem{CHQ+2006}
Canuto, C., Hussaini, M.Y., Quarteroni, A., Zang, T.A.: Spectral methods.
  Fundamentals in single domains.
\newblock Springer, Berlin (2006)

\bibitem{CPP2012}
Canuto, C., Pavarino, L.F., Pieri, A.B.: {BDDC} preconditioners for continuous
  and discontinuous {G}alerkin methods using spectral/hp elements with variable
  polynomial degree.
\newblock Submitted. (2012)

\bibitem{Nepomnyaschikh1990}
Nepomnyaschikh, S.V.: Schwarz alternating method for solving the singular
  {N}eumann problem.
\newblock Soviet J. Numer. Anal. Math. Modelling \textbf{5}(1), 69--78 (1990)

\bibitem{Nepomnyaschikh1992}
Nepomnyaschikh, S.V.: Decomposition and fictitious domains methods for elliptic
  boundary value problems.
\newblock In: D.E. Keyes et~al. (eds.) Fifth International Symposium on Domain
  Decomposition Methods for Partial Differential Equations, pp. 62--72. SIAM,
  Philadelphia (1992)

\bibitem{Oswald1996}
Oswald, P.: Preconditioners for nonconforming discretizations.
\newblock Math. Comput. \textbf{65}(215), 923--941 (1996)

\bibitem{Xu1992}
Xu, J.: Iterative methods by space decomposition and subspace correction.
\newblock SIAM Rev. \textbf{34}, 581--613 (1992)

\bibitem{Xu1996}
Xu, J.: The auxiliary space method and optimal multigrid preconditioning
  techniques for unstructured grids.
\newblock Computing \textbf{56}(3), 215--235 (1996)

\end{thebibliography}
\end{document}